\documentclass[11pt,a4paper]{article}

\usepackage{graphicx,amsmath,amssymb,enumerate,delarray}
\usepackage{macro}

\begin{document}

\begin{center}
  {\Large\bf Multiobjective optimization approach to shape and topology
  optimization of plane trusses with various aspect ratios}

\bigskip\bigskip\bigskip

{\large Makoto Ohsaki$^1*$, Saku Aoyagi$^{1}$, Kazuki Hayashi$^{1}$}

\bigskip

Department of Architecture and Architectural Engineering, Kyoto University,
Kyoto-Daigaku Katsura, Nishikyo, Kyoto 615-8540, Japan,
ohsaki@archi.kyoto-u.ac.jp

\end{center}

\bigskip\bigskip\bigskip\bigskip

\noindent{\large\bf Abstract}

  A multiobjective optimization method is proposed for obtaining the
  optimal plane trusses simultaneously for various aspect ratios of the initial
  ground structure as a set of Pareto optimal solutions generated through a
  single optimization process.
  The shape and topology are optimized simultaneously to minimize the
  compliance under constraint on the total structural volume.
  The strain energy of each member is divided into components of two coordinate
  directions on the plane.
  The force density method is used for alleviating difficulties due to existence
  of coalescent or melting nodes.
  It is shown in the numerical example that sufficiently accurate optimal
  solutions are obtained by comparison with those obtained by the linear
  weighted sum approach that requires solving a single-objective optimization
  problem many times.


\bigskip

\section{Introduction}

There have been numerous studies for topology optimization and shape
optimization of plane trusses under various static and dynamic constraints.
However, there exist several difficulties for simultaneously optimizing the
shape and topology of trusses.
The simplest approach may be to consider both of the member cross-sectional
areas and the nodal coordinates as design variables.
However, in this case, structural analysis becomes difficult when coalescent
nodes \cite{ohsaki1998} or melting nodes \cite{achtziger2007} exist, and
consequently, the stiffness matrix becomes singular.
Since topology optimization often leads to an unstable truss, 
global stability constraint is one of the important subjects for recent studies
on the shape and topology optimization, or layout optimization, of plane trusses
\cite{descamps14}.
Weldeyesus {\it et al.} \cite{weldeyesus2020} incorporated global stability
constraint as well as the bound constraints on the nodal coordinates.

The force density method (FDM) was originally developed for self-equilibrium
analysis of tension structures such as cable nets and tensegrity structures
\cite{schek1974}.
Descamps and Coelho \cite{descamps11} optimized layout of the hanging members of
a bridge truss assigning the zero lower bound for the force densities so that
all members have tension forces.
To prevent the difficulties on simultaneous shape and topology optimization of
trusses, Ohsaki and Hayashi \cite{ohsaki2017} proposed a method using the force
densities of members as design variables for the problem of minimizing the
compliance under constraint on the total structural volume.
Existence of extremely short members can be indirectly avoided by assigning
appropriate bounds on force densities \cite{shen2021}.
In the process of shape and topology optimization of plane trusses, the
locations of supports and loaded nodes are generally fixed at those specified by
the initial ground structure.
However, the optimal load path strongly depends on the ratios of distances among
the supports and loaded nodes, which are varied along with the aspect ratio of
the initial ground structure.
Therefore, it is beneficial for the engineers if the optimal trusses for various
aspect ratios of the ground structure can be found simultaneously by a single
optimization process.

Zhang and Ohsaki \cite{zhang2007} showed that the self-equilibrium state
of tensegrity structures is invariant against affine motions when
the equilibrium equations are formulated with respect to the force densities.
Utilizing this property, the force densities of members are chosen as design
variables in this study.
The total strain energy, which is equivalent to a half of the compliance, is
divided into two in-plane directional components to formulate the two objective
functions of a multiobjective optimization problem, and the optimal solutions
for various aspect ratios of the initial ground structure are generated as a set
of Pareto optimal solutions using a multiobjective genetic algorithm.

\section{FDM for simultaneous shape and topology optimization of trusses}

In this section, we briefly summarize the results in Ohsaki and Hayashi
\cite{ohsaki2017} for completeness of the paper.

Consider an optimization problem of minimizing the compliance under constraint
on the total structural volume of a truss.
Let $N_{i}$ denote the axial force of the $i$th member of a
truss consisting of $m$ members.
The truss is subjected to a single loading condition, and 
the force density $q_{i}$ is defined as the force divided by the length $L_{i}$,
i.e.,
\begin{equation}
\label{fd}
  q_{i}=\frac{N_{i}}{L_{i}} \ \ (i=1,\dots,m)
\end{equation}
Let $\bi{Q} \R^{n\times n}$ denote the force density matrix constructed from
the force densities of members and the connectivity matrix, where $n$ is the number
of nodes including the supports.
The equilibrium equations in $x$- and $y$-directions of a truss have the
following forms:
\begin{equation}
  \label{equiv}
  \begin{split}
    \bi{Q} \bi{x} = \bi{p}^{x}, \ \ 
    \bi{Q} \bi{y} = \bi{p}^{y}
  \end{split}
\end{equation}
where $\bi{x} \R^{n}$ and $\bi{y} \R^{n}$, respectively, are the vectors
of $x$- and $y$-coordinates of nodes including supports, and 
$\bi{p}^{x} \R^{n}$ and $\bi{p}^{y} \R^{n}$ are
the nodal load vectors including the reactions in $x$- and $y$-directions,
respectively, at the supports.

The nodes are classified into free and fixed nodes, which are denoted by the
subscripts `free' and `fix', respectively.
The components between the free and fixed nodes are denoted with `link'.
The set of $x$- and $y$-directional equilibrium equations in \refeq{equiv} are
partitioned as
\begin{subequations}
\label{equiv2}
  \begin{align}
  &
  \begin{pmatrix}
    \bi{Q}_{\rm free} & \bi{Q}_{\rm link} \\
    \bi{Q}^{T}_{\rm link} & \bi{Q}_{\rm fix} \\
  \end{pmatrix}
  \begin{pmatrix}
    \bi{x}_{\rm free} \\
    \bi{x}_{\rm fix} \\
  \end{pmatrix}
  =
  \begin{pmatrix}
    {\bf 0} \\
    \bi{p}^{x}_{\rm fix} \\
  \end{pmatrix} \\
  &
  \begin{pmatrix}
    \bi{Q}_{\rm free} & \bi{Q}_{\rm link} \\
    \bi{Q}^{T}_{\rm link} & \bi{Q}_{\rm fix} \\
  \end{pmatrix}
  \begin{pmatrix}
    \bi{y}_{\rm free} \\
    \bi{y}_{\rm fix} \\
  \end{pmatrix}
  =
  \begin{pmatrix}
    {\bf 0} \\
    \bi{p}^{y}_{\rm fix} \\
  \end{pmatrix}
  \end{align}
\end{subequations}
where $(\ )^{T}$ indicates the transpose of a matrix.
Note that the fixed nodes include the loaded nodes for which the locations are
fixed in the process of optimization.
Therefore, $\bi{p}^{x}_{\rm fix}$ includes the nodal loads and reactions at
supports, and no loads are applied to the free nodes.
It is seen from (\ref{equiv2}a) and (\ref{equiv2}b) that the $x$- and
$y$-coordinates of free nodes and reactions are simply scaled by $\alpha$ if the
coordinates of fixed nodes and loads are scaled by $\alpha$ in $x$- and
$y$-directions, respectively, without changing the force densities of members.
See Zhang and Ohsaki \cite{zhang2007} for details.

It is known that the absolute values of stresses of existing members of the
optimal solution take the same value, which is denoted by $\bar{\sigma}$,
because the optimal truss under single loading condition is statically
determinate and all members are fully stressed \cite{hemp1973}.
Although details are not explained here, Ohsaki and Hayashi \cite{ohsaki2017}
showed that the minimization problem of compliance under volume constraint can
be formulated as the minimization problem of the following objective function
with respect to the force density only:
\begin{equation}
  \label{obj}
  F(\bi{q}) = \ds{\sum_{i=1}^{m}
    \frac{\bar{\sigma}(L_{i}(\bi{q}))^{2} |q_{i}|}{E}}
\end{equation}
where $E$ is Young's modulus, and $L_{i}$ is a function of the force density
vector $\bi{q}$ through the nodal coordinates obtained from (\ref{equiv2}a) and
(\ref{equiv2}b).
The volume constraint is satisfied with equality by adjusting the stress level
$\bar{\sigma}$, and the cross-sectional area $A_{i}$ of member $i$ is obtained
from $L_{i}$, $q_{i}$, and $\bar{\sigma}$ using \refeq{fd} as
\begin{equation}
  A_{i} = \frac{L_{i}|q_{i}|}{\bar{\sigma}}
\end{equation}

\section{Multiobjective optimization problem}

Suppose the $i$th member is connected to two nodes $i_{1}$ and $i_{2}$.
Then, $L_{i}^{2}$ is divided into the components of $x$- and $y$-coordinates as
\begin{equation}
  \label{length}
  L_{i}^{2} = (x_{i_{2}} - x_{i_{1}})^{2} + (y_{i_{2}} - y_{i_{1}})^{2}
\end{equation}
Using this relation, the objective function \refeq{obj} is divided
into two components as 
\begin{equation}
  F(\bi{q}) = F_{x}(\bi{q}) + F_{y}(\bi{q})
\end{equation}
where
\begin{subequations}
\begin{align}
  F_{x}(\bi{q}) &= \ds{\sum_{i=1}^{m}
    \frac{\bar{\sigma}}{E}(x_{i_{2}} - x_{i_{1}})^{2}|q_{i}|} \\
  F_{y}(\bi{q}) &= \ds{\sum_{i=1}^{m}
    \frac{\bar{\sigma}}{E}(y_{i_{2}} - y_{i_{1}})^{2}|q_{i}|}
\end{align}
\end{subequations}
 
Although it is not possible to actually divide the strain energy into the
components of coordinates, $F_{x}$ and $F_{y}$ are regarded as $x$- and
$y$-components of the strain energy, respectively.
However, it should be noted that these components are different from those of
compliance which is obtained by summation of the work done by each component of
the external loads.
Therefore, for example, $F_{x}$ may have a non-zero value even when no load is
applied in $x$-direction.

Note that the axial forces computed from $\bi{q}$ and the nodal coordinates
obtained from (\ref{equiv2}a) and (\ref{equiv2}b) do not satisfy equilibrium at
the loaded nodes, because the loads are obtained as reactions at the fixed
nodes, and they are not equal to the applied loads.
In the formulation by Ohsaki and Hayashi \cite{ohsaki2017}, constraints are
given for the equilibrium of the member forces with the applied loads.
However, in this paper, the force density vector $\tilde{\bi{q}}$ satisfying
equilibrium conditions are recalculated, as follows, as a function of $\bi{q}$
at each step of optimization because it is difficult to satisfy the equality
constraints accurately when a  genetic algorithm is used for optimization in the
numerical example: 
\begin{enumerate}
  \item Assign $\bi{x}_{\rm fix}$, $\bi{y}_{\rm fix}$, and $\bi{q}$.
  \item Compute $\bi{x}_{\rm free}$ and $\bi{y}_{\rm free}$ by
    solving (\ref{equiv2}a) and (\ref{equiv2}b), respectively.
  \item Compute $L_{i}$ from the nodal coordinates and the cross-sectional area
    $A_{i}=L_{i}|q_{i}|/\bar{\sigma}$ to obtain the stiffness matrix by
    incorporating displacement boundary conditions; i.e., the displacements are
    not fixed at the loaded nodes, although locations of loaded nodes are fixed
    for optimization.
  \item Solve the stiffness equation to obtain the axial force $\tilde{N}_{i}$,
    and accordingly, the updated force density
    $\tilde{q}_{i}=\tilde{N}_{i}/A_{i}$ that is regarded as a function of
    $\bi{q}$.
\end{enumerate}

\noindent
Accordingly, $\bi{q}$ is an auxiliary variable vector that does not represent
the true force densities.
A multiobjective optimization problem is formulated, as follows, to minimize
$\tilde{F}_{x}(\bi{q})=F_{x}(\tilde{\bi{q}}(\bi{q}))$
and $\tilde{F}_{y}(\bi{q})=F_{y}(\tilde{\bi{q}}(\bi{q}))$
under the bound constraints of $\bi{q}$:
\begin{subequations}
\label{p2}
\begin{alignat}{2}
 & \mathrm{Minimize} & \ \
 & \tilde{F}_{x}(\bi{q}) \ {\rm and} \ \tilde{F}_{y}(\bi{q}) \\
 & \mathrm{subject\,to} &
 & q_{i}^{\rm L} \le q_{i} \le q_{i}^{\rm U} \ \ (i=1,\dots,m)
\end{alignat}
\end{subequations}
where $q_{i}^{\rm U}$ and $q_{i}^{\rm L}$ are the upper and lower bounds of
$q_{i}$, respectively.

In the following numerical examples, Pareto optimal solutions of Problem
\refeq{p2} are generated simultaneously through a single optimization process
using a multiobjective genetic algorithm.
The problem is also solved using the linear weighted sum approach for comparison
purpose.
The objective function $F^{*}(\bi{\bi{q}})$ is formulated using nonnegative
weight coefficients $\mu_{x}$ and $\mu_{y}$ as
\begin{equation}
\label{ff}
\begin{split}
  F^{*}(\bi{\bi{q}})
  &= \mu_{x}\tilde{F}_{x}(\bi{q}) + \mu_{y}\tilde{F}_{y}(\bi{q}) \\
  &= \ds{\sum_{i=1}^{m} \frac{\bar{\sigma}}{E}}
    [\mu_{x}(x_{i_{2}} - x_{i_{1}})^{2}
    + \mu_{y}(y_{i_{2}} - y_{i_{1}})^{2}]\sqrt{q_{i}^{2}+c}
\end{split}
\end{equation}
where $c$ is a small positive value for smoothing the function of absolute value
for the optimization process using a nonlinear programming approach. 

By comparing $F^{*}(\bi{\bi{q}})$ in \refeq{ff} and $F(\bi{q})$ in \refeq{obj}
with the relation \refeq{length}, we can see that the optimal solution of a
truss after scaling the locations of the supports and the loaded nodes in the
initial ground structure by the factors $\sqrt{\mu_{x}}$ and $\sqrt{\mu_{y}}$,
respectively, in $x$- and $y$-directions is found by minimizing
$F^{*}(\bi{\bi{q}})=\mu_{x}\tilde{F}_{x}(\bi{q})+\mu_{y}\tilde{F}_{y}(\bi{q})$. 
Thus, the optimal solutions of the trusses with various aspect ratios can be
found by solving a multiobjective optimization problem of minimizing 
$\tilde{F}_{x}(\bi{q})$ and $\tilde{F}_{y}(\bi{q})$ only once using an
appropriate method, e.g., multiobjective genetic algorithm that can generate the
set of Pareto optimal solutions with single process of optimization.

\section{Numerical examples}

\begin{figure}[tb]
  \centering
  \includegraphics[width=0.35\textwidth]{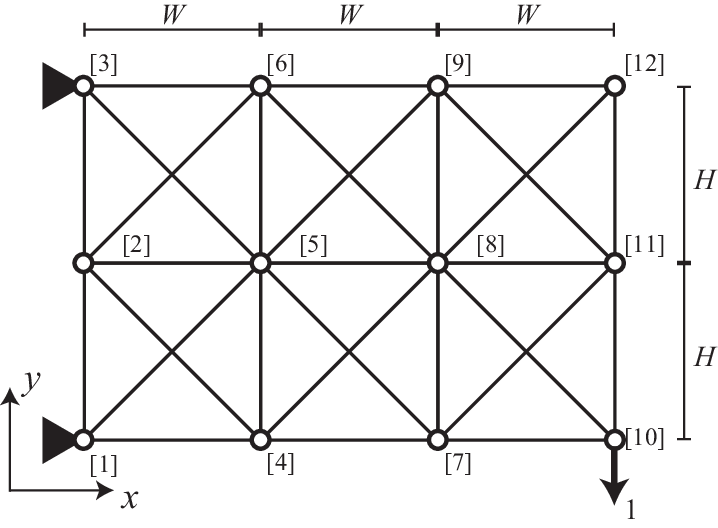}
  \caption{$3\times 2$ grid truss.}
  \label{model32}
\end{figure}

Optimal shapes and topologies are found for a $3\times 2$ grid truss as shown
in \reffig{model32} as the initial ground structure, where a unit load is
applied in negative $y$-direction at node 10.
Note that the crossing diagonal members are not connected at their centers, and
accordingly, the total number of members is 29.
The scale in $x$- and $y$-directions are varied with the scales $\sqrt{\mu_{x}}$
and $\sqrt{\mu_{y}}$, respectively, i.e., $H$ and $W$ in \reffig{model32} are 
scaled to $\mu_{x}W$ and $\mu_{y}H$, respectively.
In the following, ratio $r=\sqrt{\mu_{x}/\mu_{y}}$ is taken as a parameter.

Pareto optimal solutions are found using the non-dominated sorting genetic
algorithm (NSGA-III) implemented in the Python library
DEAP Ver. 1.3.1 \cite{DEAP}.
The numbers of population and generations are 40 and 500, respectively,
and the probabilities of crossover and mutation are 0.9 and 0.0345 $(=1/m)$,
respectively.
Young's modulus is $E=1.0$.
The stress of the existing members is tentatively set as $\bar{\sigma}=1.0$, and
the cross-sectional areas are uniformly scaled to have the specified total
structural volume of 100 after obtaining the Pareto optimal solutions.
Note that the units are omitted because they are not important.

Each individual is coded using real numbers $q_{i}$ $(i=1,\dots,m)$ as variables.
Let $q_{i0}$ denote the force density of the $i$th member of the ground
structure in \reffig{model32} with uniform cross-section.
In the initial population, individuals are generated randomly in the range
\begin{equation}
  q_{i0}-10\le q_{i}\le q_{i0}+10
\end{equation}
After generating the Pareto front obtained using NSGA-III, the weight
coefficients for each solution are approximately computed using the slope of the
tangent line of the Pareto front in the objective function space.
Let $\bi{S}_{i}=(F_{xi},F_{yi})$ denote the location of the $i$th Pareto optimal
solution aligned in the increasing order of $F_{x}$ on the $(F_{x},F_{y})$
plane.
The slope of the line connecting $\bi{S}_{i}$ and $\bi{S}_{i-1}$ is denoted by
$\beta_{i}$, i.e.,
\begin{equation}
  \beta_{i} = \frac{F_{yi}-F_{y(i-1)}}{F_{xi}-F_{x(i-1)}}
\end{equation}

Then the ratio between the weight coefficients $\mu_{xi}$ and $\mu_{yi}$ for the
$i$th solution is chosen from the range between $\beta_{i+1}$ and
$\beta_{i}$ to estimate the value $r_{i}=\sqrt{\mu_{xi}/\mu_{yi}}$.

To compare the solutions to those obtained by NSGA-III, the sequential quadratic
programming in the library SNOPT Ver. 7 \cite{gill1997} is used for minimizing the
weighted sum of the objective functions $F^{*}(\bi{\bi{q}})$ in \refeq{ff} with
the smoothing parameter $c=1.0\times 10^{-10}$.
Both for the genetic algorithm and linear weighted sum approach, the lower and
upper bounds for $q_{i}$ are given as $q_{i}^{\rm L}=-1000$ and
$q_{i}^{\rm U}=1000$, respectively, which have sufficiently large absolute
values not to restrict the solution space.
Three individuals are generated by assigning three different sets of
weight coefficients corresponding to $r=1.0$, 2.0, and 3.0 for the linear
weighted sum approach, and they are added to the initial population of NSGA-III
to enhance diversity of the Pareto optimal solutions.

\begin{figure}[tb]
  \centering
  \includegraphics[width=0.5\textwidth]{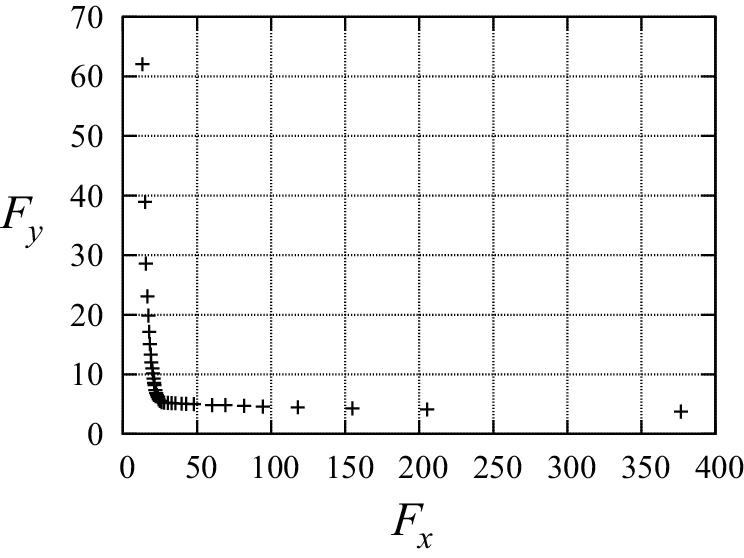}
  \caption{Pareto optimal solutions in the objective function space
    $(F_{x},F_{y})$.}
  \label{pareto-1}
\end{figure}


\begin{figure}[tb]
  \centering
  \includegraphics[width=0.5\textwidth]{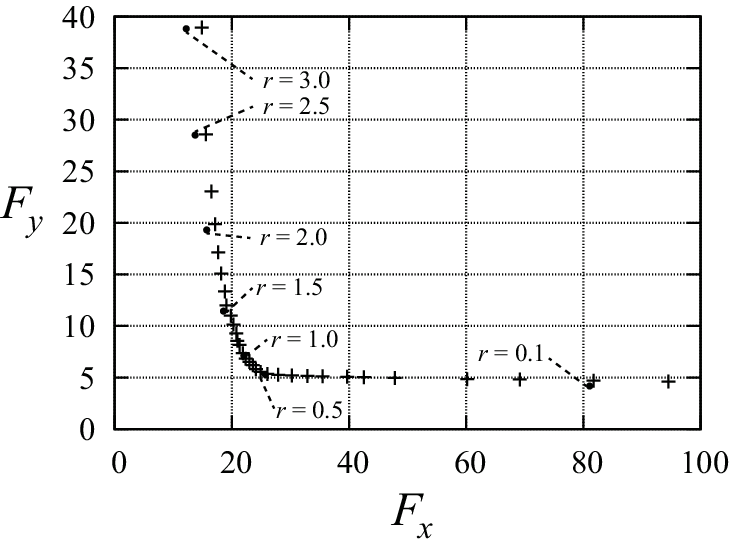}
  \caption{Pareto optimal solutions in the objective function space
    $(F_{x},F_{y})$; `+': solutions by NSGA-III,
    dots: solutions of Problem \refeq{p2}.}
  \label{pareto-3}
\end{figure}

\begin{figure}[tb]
\begin{center}
\begin{minipage}{0.15\textwidth}
  \centering
  \includegraphics[scale=0.3]{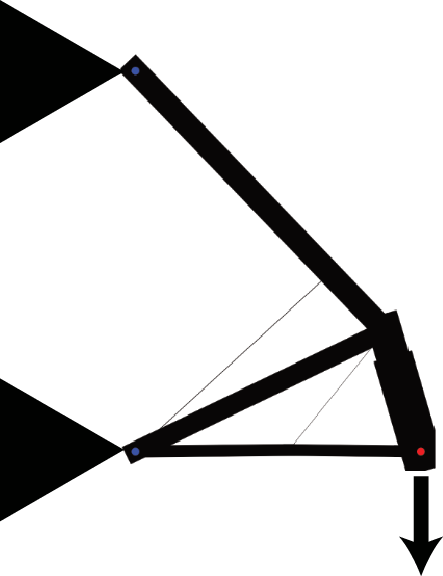}
  \par (a) \par
\end{minipage}
\hspace{2cm}
\begin{minipage}{0.3\textwidth}
  \centering
  \includegraphics[scale=0.3]{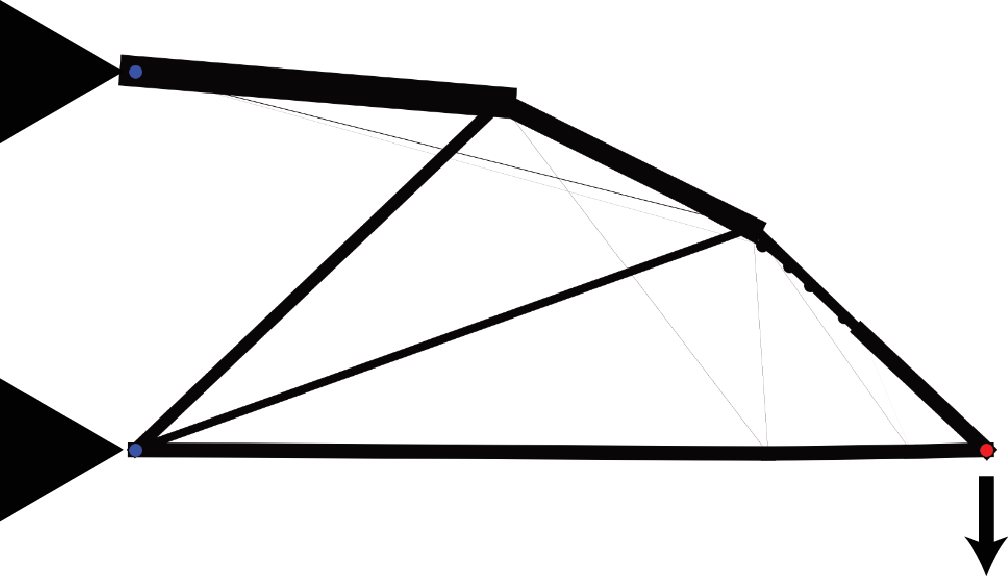}
  \par (b) \par
\end{minipage}
\end{center}
\par
\begin{center}
\begin{minipage}{0.44\textwidth}
  \centering
  \includegraphics[scale=0.3]{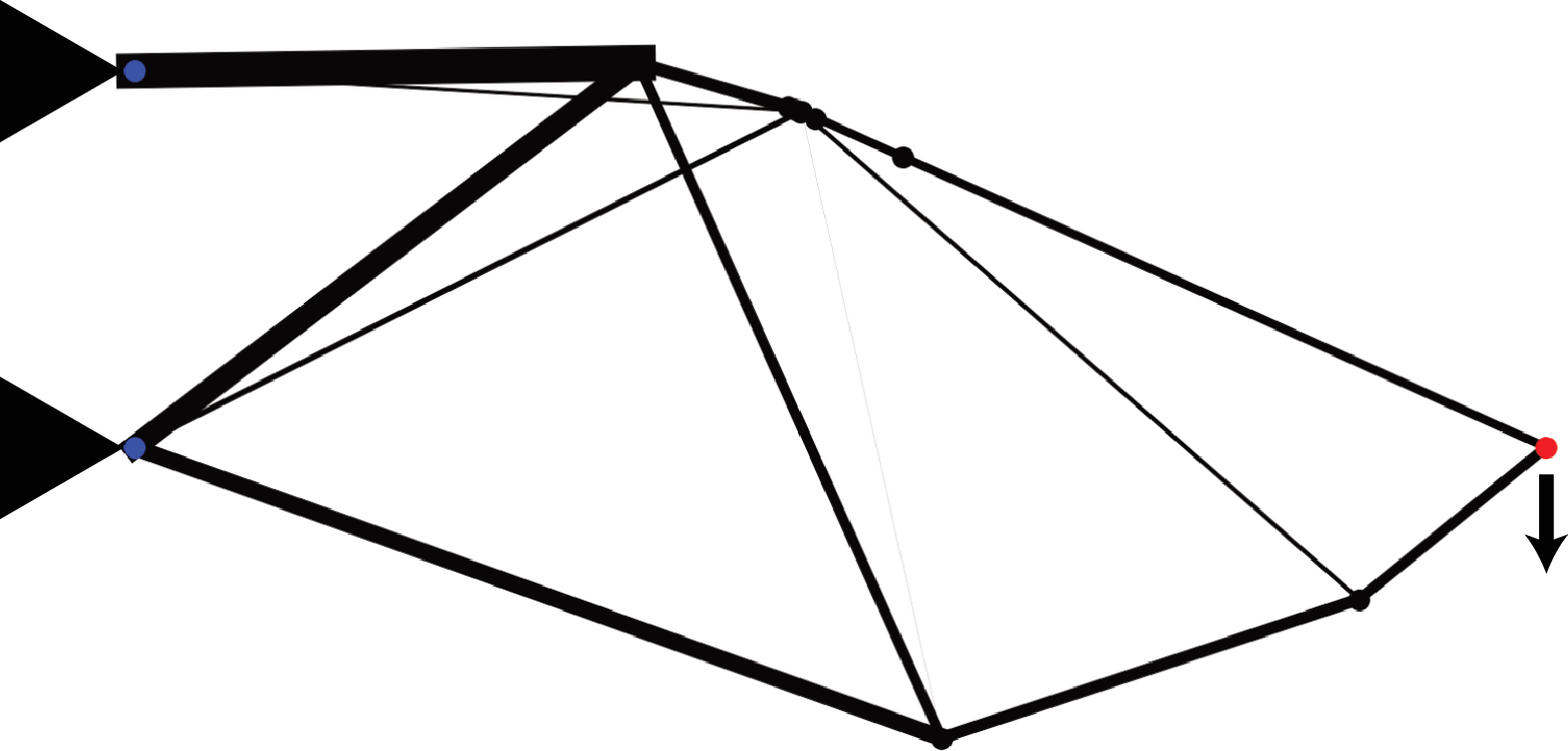}
  \par (c) \par
\end{minipage}
\end{center}
\caption{Optimal solutions: (a)$r=0.5$, (b) $r=1.5$, (c) $r=2.5$.}
\label{opt}
\end{figure}

The `+' marks in \reffig{pareto-1} are the Pareto optimal solutions obtained by
NSGA-III in the objective function space, where the solutions outside of the
ranges of this figure are excluded.
Note that optimization using NSGA-III has been carried out ten times from 
different initial solutions.
The numbers of Pareto optimal solutions are 40, which is equal to the number of
population, for all cases.
\reffig{pareto-1} shows the best case in view of diversity of solutions.
It is seen from \reffig{pareto-1} that sufficiently large number of solutions
are obtained, and the solutions are located along a convex curve in the
objective function space.
Figure \ref{pareto-3} shows the Pareto optimal solutions in a smaller range of
objective function space, where the optimal solutions obtained by the weighted
sum approach are also plotted with dots and the corresponding $r$ values.
It is seen from the figure that the objective functions obtained by NSGA-III are
slightly larger for a large value of $r$; however, a good
agreement in the objective values has been observed.
The solutions obtained by NSGA-III for $r=0.5$, 1.5, 2.5 are shown in
Figs.~\ref{opt}(a), (b), (c), respectively, with the same
span in $y$-direction after scaling
the nodal coordinates by 0.5, 1.5, and 2.5 in $x$-direction.
Note that the width of each member is proportional to its cross-sectional area
after scaling the ground structure to have the total structural volume of 100.

\begin{table}[tb]
  \centering
  \caption{Optimal compliance values obtained by three different methods.}
  \begin{tabular}{cccc} \hline
    Method&$r=0.5$&$r=1.5$&$r=2.5$\\\hline
    Scaling&1.38&31.53&155.10\\
    Weighted sum&1.46&31.78&151.29\\
    NSGA-III&1.39&31.38&158.28\\\hline
  \end{tabular}
  \label{obj-val}
\end{table}

\reftab{obj-val} shows the optimal compliance values obtained by the three
different methods, where `Scaling' refers to the solution 
obtained after actually scaling in $x$-direction by the
ratio $r$ and solving the single-objective optimization problem using the method
by Ohsaki and Hayashi \cite{ohsaki2017}.
A good agreement is confirmed in \reftab{obj-val}, although inaccuracy in the
$r$ value for NSGA-III estimated from the slope of Pareto front is also regarded as a source of difference among the
compliance values obtained by the three methods. 

\section{Conclusions}

A new method of simultaneous optimization of geometry and topology has been 
presented for plane trusses to minimize the compliance under constraint on
the total structural volume.
Using the proposed method, the optimal shape and topology of trusses for various
aspect ratios of the initial ground structure can be obtained as a set of Pareto
optimal solutions through a single optimization process of a multiobjective
genetic algorithm.
The ratio between the weight coefficients for the two components of strain
energy is obtained from the slope of the objective function values of the Pareto
front in the objective function space.
The ratio of scaling factors in the two directions of the ground structure is
computed as the square root of the ratio of weight coefficients.
It has been confirmed in the numerical example that the solutions corresponding
to various aspect ratios of the initial ground structure are obtained with good
accuracy in comparison to the solutions obtained by the weighted sum approach
that demands solving an optimization problem for each aspect ratio of the
ground structure.

\end{document}